\let\und=\underline

\def\MM{{\overline{M}}}

\def\M0{{\overline{M}_{0,n}(\mathbb{P}^{r},d)}}

\def\mm{{\overline{\mathcal{M}}}}
\def\m0{{\overline{\mathcal{M}}_{0,n}(\mathbb{P}^{r},d)}}

\def\mmm{{\overline{\mathcal{M}}_{g,n}(X,\beta)}}

\def\nobarm0{{\mathcal{M}_{0,n}(\mathbb{P}^{r},d)}}

\def\D{{\Delta}}
\def\G{{\Gamma}}

\def\I{{\mathcal I}}

\def\O{{\mathcal O}}

\def\AA{{\mathbb A}}
\def\CC{{\mathbb C}}

\def\NN{{\mathbb N}}
\def\PP{{\mathbb P}}
\def\QQ{{\mathbb Q}}
\def\R{{\mathcal R}}

\def\C{{\mathcal C}}

\def\a{\alpha}
\def\b{\beta}

\def\la{\lambda}

\def\g{\gamma}

\def\r{\rho}
\def\s{\sigma}
\def\t{\tau}
\def\w{\omega}

\def\pr{^{\prime}}

\def\x{\times}
\def\*{\otimes}
\def\iso{\simeq}
\def\sub{\subset}

\def\+{\oplus}

\def\bra{\langle}
\def\ket{\rangle}

\def\ra{\rightarrow}

\def\pj{\operatorname{pr}}
 
\def\ker{\operatorname{ker}}

\def\ev{\operatorname{ev}}
\def\dim{\operatorname{dim}}

\def\deg{\operatorname{deg}}

\def\Bl{\operatorname{B\ell}}
\def\eqeul{\operatorname{Euler_T}}

\def\val{\operatorname{val}}

\def\Aut{\operatorname{Aut}}

\documentclass[12pt]{article}

\usepackage{amssymb}
\usepackage{latexsym,amsfonts}
\usepackage{amsthm}
\usepackage{pstricks,pst-node,pst-tree}
\usepackage{graphicx}
\usepackage{amstext}
\usepackage{amsmath}

\begin{document}

\newtheorem{thm}{Theorem}
\newtheorem{cor}{Corollary}
\newtheorem{Def}{Definition}
\newtheorem{eg}{Example}
\newtheorem{prop}{Proposition}
\newtheorem{rmk}{Remark}
\newtheorem{lem}{Lemma}
\newtheorem{conj}{Conjecture}

\title{
A presentation for the Chow ring of
$\mm_{0,2}(\PP^1,2)$}

\author{Jonathan Cox
}

\date{\today}

\maketitle

\section{Introduction}
\label{sec:intro}

In a previous article (\cite{C5}), we computed the Betti numbers and an additive basis for the Chow ring
$A^*(\mm_{0,2}(\PP^r,2))$ of the moduli space of degree two stable maps from 2-pointed curves to projective
space of arbitrary dimension $r$. These are the first steps in a program to give presentations for these Chow 
rings. In the present paper, we will complete this program for the special case $r=1$ by showing that
\[A^*(\mm_{0,2}(\PP^1,2))\iso\frac{\QQ[D_0,D_1,D_2,H_1,H_2,\psi_1,\psi_2]}
{\left(\begin{array}{c}
H_1^2, H_2^2,D_0\psi_1,D_0\psi_2,D_2-\psi_1-\psi_2, \\
\psi_1-\frac{1}{4}D_1-\frac{1}{4}D_2-D_0+H_1, (D_1+D_2)^3, \\
\psi_2-\frac{1}{4}D_1-\frac{1}{4}D_2-D_0+H_2, D_1\psi_1\psi_2
\end{array}\right)}
\]
This is the first known presentation for a Chow ring of a 
moduli space of stable maps of degree greater than one with more than one 
marked point. 
 
We begin in Section \ref{sec:gen} by listing some natural divisor classes that occur
in the Chow rings of all moduli spaces of stable maps
to projective spaces. Although these divisor classes  
do not generate such Chow rings in general,
we will demonstrate that they do in the case of $\mm_{0,2}(\PP^1,2)$.
Relations among these classes are found in Section \ref{sec:rel}.
In the course of proving these relations, computations on moduli spaces
of stable maps with fewer marked points or lower degree naturally arise.
Section \ref{sec:simpler} describes presentations for the Chow rings of these
simpler moduli spaces.  All of these pieces are compiled to give
the whole presentation in Section \ref{sec:prez}. This Section also argues
that the presentation is complete, {\em ie.}, it doesn't leave out any essential
generators.
Finally, Section \ref{sec:app} applies the presentation to
compute the genus zero, degree two, two-pointed gravitational
correlators of $\PP^1$. Algorithms for computing theses values have previously
been developed; see \cite{KM2} and \cite{CK}, for example. However, the method given here is more explicit.

We will work over the field $\CC$ of complex numbers. As in \cite{C5}, we
let $\und{n}=\NN\cap[1,n]$ be the initial segment consisting of the first $n$ natural numbers.
We will also continue to freely make use of the homology isomorphism
$A^*(\m0)\ra H^*(\m0)$.

\subsection{Acknowledgements}

The content of this article derives from a part of my doctoral dissertation.
I am deeply grateful to my dissertation adviser, Dr. Sheldon Katz, for 
presenting me with many opportunities, for plenteous financial support, for insight, 
for encouragement, and for inspiration. 
I am further beholden to Dr. William Jaco and Dr. 
Alan Adolphson for providing 
additional financial 
support during work on this project.
I appreciate the generous support given to me by the University of Illinois 
mathematics department during my two and a half years 
as a visiting graduate student there. I also acknowledge 
with humble gratitude the
Oklahoma State University mathematics department for extended
support during that time. 

\section{Generators}
\label{sec:gen}
This section describes three types of divisor classes found in all
Chow rings \linebreak[4] $A^*(\m0)$. We will
see in Section \ref{sec:comp} that divisors 
of these types generate the Chow
ring of the moduli space $\mm_{0,2}(\PP^1,2)$. All of these divisors
have equivariant lifts, produced by taking the equivariant first Chern classes
of their corresponding equivariant line bundles. We will use the same
symbols for these equivariant versions; the meaning should be clear from
the context.

First we will describe the {\em hyperplane pullbacks}.
Let $\ev_1,\ldots,\ev_n$ be the evaluation maps on $\mmm$, where $\ev_i$ takes a
stable map $(C,x_1,\ldots,x_n,f)$ to $f(x_i)$. 
Any cohomology class on $X$ pulls back under any evaluation
map to a class on $\mmm$, which
may be considered as an element in the Chow ring under homology isomorphism.
In the case $X=\PP^r$, we 
can pull the hyperplane class $H$ back under each evaluation, 
getting the $n$ hyperplane pullbacks $H_i=\ev_i^*(H)$.  

Next we will describe the {\em boundary divisors}.
The boundary of $\mm_{0,n}(X,\b)$ by definition consists
of the locus of stable maps with reducible domain curves. 
It is a divisor with normal crossings, and its irreducible components
are in 1--1 correspondence with quadruples $(A,d_A,B,d_B)$,
where $A,B\sub\und{n}$ partition $\und{n}$, $d_1+d_2=d$, and if 
$d_A=0$ (resp. $d_B=0$), then $A$ (resp. $B$) has at least two elements. 
Such a boundary divisor and its
class in the Chow ring are both denoted $D_{A,d_A,B,d_B}$. Geometrically,
the divisor $D_{A,d_A,B,d_B}$ corresponds to the closure of the
locus of stable maps where the domain curve has two 
components, one having marked points labeled by $A$ and mapping to $\PP^r$
with degree $d_A$, and the other having marked points labeled by $B$ and 
mapping to $\PP^r$ with degree $d_B$. We represent this divisor by the
following picture.

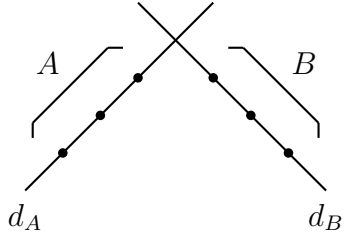
\begin{center}
\begin{pspicture}(0,0)(5,3.5)
\pnode(0.5,.5){a}
\dotnode(1,1){b}
\dotnode(1.5,1.5){c}
\dotnode(2,2){d}
\pnode(3,3){e}
\pnode(2,3){f}
\dotnode(3,2){g}
\dotnode(3.5,1.5){h}
\dotnode(4,1){i}
\pnode(4.5,0.5){j}
\pnode(.6,1.2){k}
\pnode(.6,1.4){l}
\pnode(1.6,2.4){m}
\pnode(1.8,2.4){n}
\pnode(4.4,1.2){o}
\pnode(4.4,1.4){p}
\pnode(3.4,2.4){q}
\pnode(3.2,2.4){r}
\ncline{a}{e}
\ncline{f}{j}
\ncline{k}{l}
\ncline{l}{m}\naput{$A$}
\ncline{m}{n}
\ncline{o}{p}
\ncline{q}{p}\naput{$B$}
\ncline{q}{r}
\uput{5pt}[d](.5,.5){$d_A$}
\uput{5pt}[d](4.5,.5){$d_B$}
\end{pspicture}
\end{center}
Note that the domain curves of stable maps lying in a boundary
divisor may have more than two components. In the limit, some combinations
of marked points and the node may coincide, causing new components to sprout.
The marked points involved in the ``collision" will appear on this new
component. Observe that, although marked points can thus
migrate to newly added
components in the limit, they cannot move onto components that already existed
as long as the node separating the components is maintained.
Additionally, the map itself may degenerate in such a way that the number of
components of the domain curve increases. 
The diagram representation given
above for divisors can easily be extended to describe the closures of
other degeneration loci. This description is an alternative to the
dual graphs used to describe the degeneration loci in \cite{C5}. We will use  the diagram
representation above when referring to their closures of the degeneration loci as opposed
to the loci themselves. Such a diagram therefore
directly describes only a generic element of the locus it represents. For
example, the diagram

\begin{center}
\begin{pspicture}(0,0)(5,2)
\pnode(1,1){a}
\dotnode(2,1){c}
\dotnode(3,1){d}
\pnode(4,1){b}
\ncline{a}{b}
\uput{5pt}[r](4,1){2}
\end{pspicture}
\end{center}
represents (a generic element of) $\mm_{0,2}(\PP^r,2)$ itself.

In $A^*(\mm_{0,2}(\PP^r,2))$, there are exactly three boundary divisors.
We use the notation  
$D_0=D_{\underline{2},0,\emptyset,2}$, $D_1=D_{\underline{2},1,\emptyset,1}$, 
and
$D_2=D_{\underline{1},1,\{2\},1}$ for these divisors.  
Thus the domain of a generic stable map in $D_0$ has one collapsed component 
containing both marked points and one component of degree two. 
Generic elements of $D_1$ are maps which have degree one on each of the two
components of the domain curve, with both 
marked points lying on the same component.
Finally, $D_2$ is the boundary
divisor whose generic maps have degree one on each component with one marked 
point on each component. Note that $D_i$ corresponds to curves with $i$
degree one components containing marked points; this aids in remembering the 
notation.
The diagrams for these three divisors 
follow.

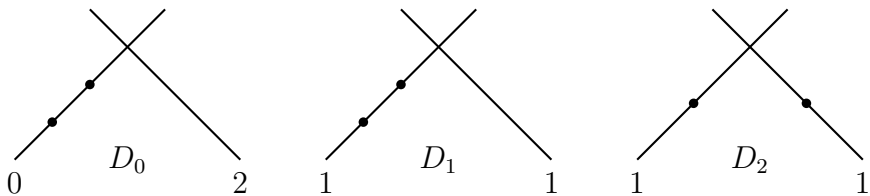
\begin{center}
\begin{pspicture}(0,0)(4,3)
\rput(2,0.5){$D_0$}
\pnode(0.5,.5){a}
\dotnode(1.5,1.5){c}
\dotnode(1,1){d}
\pnode(2.5,2.5){e}
\pnode(1.5,2.5){f}
\pnode(3.5,0.5){j}
\ncline{a}{e}
\ncline{f}{j}
\uput{5pt}[d](.5,.5){0}
\uput{5pt}[d](3.5,.5){2}
\end{pspicture}
%d1
\begin{pspicture}(0,0)(4,3)
\rput(2,0.5){$D_1$}
\pnode(0.5,.5){a}
\dotnode(1.5,1.5){c}
\dotnode(1,1){d}
\pnode(2.5,2.5){e}
\pnode(1.5,2.5){f}
\pnode(3.5,0.5){j}
\ncline{a}{e}
\ncline{f}{j}
\uput{5pt}[d](.5,.5){1}
\uput{5pt}[d](3.5,.5){1}
\end{pspicture}
%d2
\begin{pspicture}(0,0)(4,3)
\rput(2,0.5){$D_2$}
\pnode(0.5,.5){a}
\dotnode(1.25,1.25){c}
\dotnode(2.75,1.25){d}
\pnode(2.5,2.5){e}
\pnode(1.5,2.5){f}
\pnode(3.5,0.5){j}
\ncline{a}{e}
\ncline{f}{j}
\uput{5pt}[d](.5,.5){1}
\uput{5pt}[d](3.5,.5){1}
\end{pspicture}
\end{center}

We use $D_1$ as an example to 
illustrate the further degeneration that can occur within a boundary
divisor. Contained within $D_1$ are loci with the following diagrams.

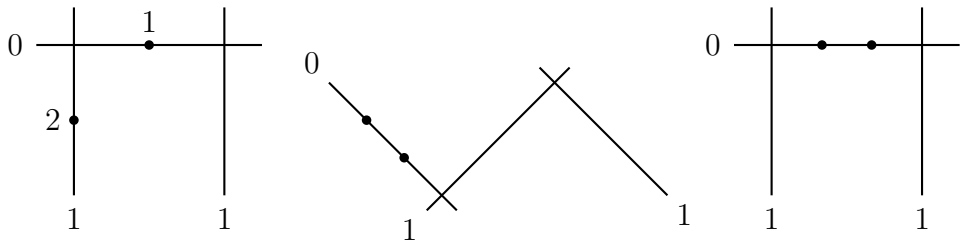
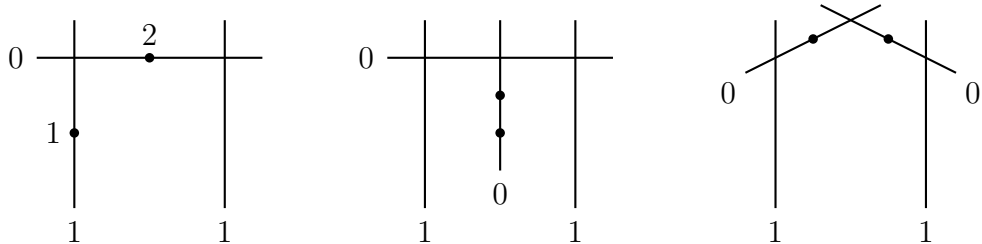
\begin{center}
%101(1)
\begin{pspicture}(0,0)(3.75,3.5)
\pnode(1,.5){a}
\pnode(1,3){b}
\pnode(.5,2.5){c}
\pnode(3.5,2.5){d}
\pnode(3,3){e}
\pnode(3,.5){f}
\dotnode(1,1.5){i}
\dotnode(2,2.5){j}
\ncline{a}{b}
\ncline{c}{d}
\ncline{e}{f}
\uput{5pt}[d](1,.5){1}
\uput{5pt}[l](.5,2.5){0}
\uput{5pt}[d](3,.5){1}
\uput{5pt}[l](1,1.5){2}
\uput{5pt}[u](2,2.5){1}
\end{pspicture}
%011
\begin{pspicture}(0,0)(5.25,2.5)
\pnode(.5,2){a}
\pnode(2.2,.3){b}
\pnode(1.8,.3){c}
\pnode(3.7,2.2){d}
\pnode(3.3,2.2){e}
\pnode(5,.5){f}
\dotnode(1,1.5){i}
\dotnode(1.5,1){j}
\ncline{a}{b}
\ncline{c}{d}
\ncline{e}{f}
\uput{5pt}[ul](.5,2){0}
\uput{5pt}[dl](1.8,.3){1}
\uput{5pt}[dr](5,.5){1}
\end{pspicture}
%101(b)
\begin{pspicture}(0,0)(3.75,3.5)
\pnode(1,.5){a}
\pnode(1,3){b}
\pnode(.5,2.5){c}
\pnode(3.5,2.5){d}
\pnode(3,3){e}
\pnode(3,.5){f}
\dotnode(1.67,2.5){i}
\dotnode(2.33,2.5){j}
\ncline{a}{b}
\ncline{c}{d}
\ncline{e}{f}
\uput{5pt}[d](1,.5){1}
\uput{5pt}[l](.5,2.5){0}
\uput{5pt}[d](3,.5){1}
\end{pspicture}

%101(2)
\begin{pspicture}(0,0)(3.75,3.5)
\pnode(1,.5){a}
\pnode(1,3){b}
\pnode(.5,2.5){c}
\pnode(3.5,2.5){d}
\pnode(3,3){e}
\pnode(3,.5){f}
\dotnode(1,1.5){i}
\dotnode(2,2.5){j}
\ncline{a}{b}
\ncline{c}{d}
\ncline{e}{f}
\uput{5pt}[d](1,.5){1}
\uput{5pt}[l](.5,2.5){0}
\uput{5pt}[d](3,.5){1}
\uput{5pt}[l](1,1.5){1}
\uput{5pt}[u](2,2.5){2}
\end{pspicture}
\hspace{.25in}
%1001tail
\begin{pspicture}(0,0)(3.75,3.5)
\pnode(1,.5){a}
\pnode(1,3){b}
\pnode(.5,2.5){c}
\pnode(3.5,2.5){d}
\pnode(3,3){e}
\pnode(3,.5){f}
\pnode(2,1){g}
\pnode(2,3){h}
\dotnode(2,1.5){i}
\dotnode(2,2){j}
\ncline{a}{b}
\ncline{c}{d}
\ncline{e}{f}
\ncline{g}{h}
\uput{5pt}[d](1,.5){1}
\uput{5pt}[l](.5,2.5){0}
\uput{5pt}[d](3,.5){1}
\uput{5pt}[d](2,1){0}
\end{pspicture}
\hspace{.25in}
%1001
\begin{pspicture}(0,0)(3.75,4.5)
\pnode(1,.5){a}
\pnode(1,3){b}
\pnode(.6,2.3){c}
\pnode(2.4,3.2){d}
\pnode(1.6,3.2){e}
\pnode(3.4,2.3){f}
\pnode(3,.5){g}
\pnode(3,3){h}
\dotnode(1.5,2.75){i}
\dotnode(2.5,2.75){j}
\ncline{a}{b}
\ncline{c}{d}
\ncline{e}{f}
\ncline{g}{h}
\uput{5pt}[d](1,.5){1}
\uput{5pt}[dl](.6,2.3){0}
\uput{5pt}[dr](3.4,2.3){0}
\uput{5pt}[d](3,.5){1}
\end{pspicture}
\end{center}

\noindent
The marked points are not labeled in diagrams where the distinction does
not affect the boundary class, either because both marked points lie on
the same component or because of symmetry.

These three boundary divisor classes 
together with the hyperplane pullbacks $H_1$ and $H_2$ generate the linear
part of the ring $A^*(\mm_{0,2}(\PP^1,2))$. We can see this by using the additive
basis from \cite{C4}. Recall that there the moduli space $\mm_{0,2}(\PP^r,2)$ was 
stratified as

\[\mm_{0,2}(\PP^1,2)=U\coprod V_1\coprod V_2\coprod Y\text{,}\]
where in the case $r=1$ we have $U\iso(\PP^1)^2\times\AA^2$, $V_i\iso(\PP^1)^2\times\AA^1$, and
$Y\iso\PP^1\times[\PP^1/S_2]$. Now $U$ is the locus where the domain curve has an irreducible
component mapped with degree two. Since its $(\PP^1)^2$ factor parametrizes the images of the
marked points, it is easy to see that the two independent ring generators of $A^*(U)$ correspond to $H_1$ and $H_2$
in the notation above. Since $V_1$ densely contains the locus where both marked points are on a
degree one component and $V_2$ densely contains the locus where each degree one component has a marked
point, we can conclude that $V_i$ corresponds to $D_i$ in the notation above. Since the relevant linear classes
are independent in the additive basis and the first Betti number of $\mm_{0,2}(\PP^1,2)$ is four, we have demonstrated
that even the classes $H_1$, $H_2$, $D_1$, and $D_2$ generate its Picard group.
We will give an additional argument in Section \ref{sec:lla} that {\em any} four of these five divisor classes 
suffice to generate.
Even better, we will see 
in Section \ref{sec:comp} that these classes
generate this entire Chow ring.

The last collection of divisor classes consists of the {\em $\psi$-classes}.
Let $\pi:\mathcal{C}\ra\mm$ be the universal curve of a moduli stack $\mm=\m0$, and let 
$\s_1,\ldots,\s_n$ be the $n$ universal sections.
We define the cotangent line bundles $L_1,\ldots,L_n$ by $L_i=s_i^*(\w_{\pi})$, where
$\w_{\pi}$ is the relative dualizing sheaf of the universal curve.
Since the total space is a smooth stack, we can write
\[\w_\pi=K_\mathcal{C}\*\pi^*K_{\mm}^\vee\text{,}\]
where $K_\mathcal{C}$ and $K_{\mm}$ are the canonical bundles.
At a point $b\in \mm$, the fiber of $L_i$ is the cotangent space to the curve
$\mathcal{C}_b$ at the point $s_i(b)$. We also define the {\em $\psi$-class}
$\psi_i$ to be the first Chern class $c_1(L_i)$ for each $i$.
It is straightforward to check that these $\psi$-classes are
universal as well: given any morphism $g:S\ra\mm$, the pullbacks 
$g^*(\psi_1),\ldots,g^*(\psi_n)$ are the $\psi$-classes on the induced family.

Although the $\psi$-classes are not strictly necessary as generators for
the ring $A^*(\mm_{0,2}(\PP^1,2))$, we include them because their geometric
nature makes some of the relations much easier to understand and state.
We will give here one example of the usefulness of the $\psi$-classes
in describing geometric conditions.

To make this example easier to state, we first introduce a slight
modification of the concept of $\psi$-classes.
Restricting to the closure of a particular degeneration locus, let $p$
be a node at the intersection of components $E_1$ and $E_2$, as in the
diagrams earlier in this section. Then we can define classes
$\tilde{\psi}_{p,E_i}$ essentially just like we defined $\psi$-classes.  The
only difference here is that we
additionally specify which branch to consider $p$ to be lying on, 
so that the cotangent
space is one-dimensional. In Theorem \ref{norm}, this type of class will be
denoted by $e_F$, where $F$ is the corresponding flag in the graph of the stable map. 
If, we remove, say, $E_2$ and then the 
associated connected component of the curve, replacing the node with an
auxiliary marked point $s_{\bullet}$, then $\tilde{\psi}_{p,E_1}$ becomes
a legitimate $\psi$-class $\psi_{\bullet}$ on the resulting moduli space.

An important fact which will be used here many times 
is that a collapsed rational component with exactly three special points is a
{\em rigid object}. In other words, in deformations of a stable map, the marked points and
nodes on such a component
cannot be moved around internally on the component. 
Among other things, this says that, once such a component appears in
a degeneration, it will remain in any further degeneration. We have already
used this fact implicitly, for example, in describing the possible further 
degenerations of $D_1$ above. This rigidity is well-known and is
intuitively clear from a brief study of the automorphism group of $\PP^1$. 
Morover, the rigidity of such a component is equivalent to the 
vanishing of the corresponding $\psi$-classes and $\tilde{\psi}$-classes
on that component.
For $\psi$-classes, this is because the cotangent 
line bundles are trivial if and only if the marked points are fixed.
A similar statement holds for $\tilde{\psi}$-classes. 
For us, the most important special case is the vanishing of
$\psi_1$ and $\psi_2$ on $D_0$, which will be described in detail in Section \ref{sec:geomrel}.

\section{Presentations for the Chow rings of some simpler spaces}
\label{sec:simpler}
%\label{sec:0n11}

The study of degree one
stable maps to $\PP^1$ is relatively simple because every degree one
morphism from $\PP^1$ to $\PP^1$ is an isomorphism. It follows that 
all degree one stable 
maps from $\PP^1$ to itself are isomorphic to the 
identity. More generally, genus 
zero, degree one stable maps never have nontrivial automorphisms. Thus the
corresponding moduli spaces may just as well be considered as fine moduli
schemes, as no loss of information results.
First, we mention that
the moduli space $\mm_{0,0}(\PP^1,1)$ is a point because the domain
of such a stable map is always $\PP^1$.
Second, note that 
$\mm_{0,3}(\PP^1,0)\iso\MM_{0,3}\x\PP^1\iso\PP^1$ since
the moduli space of stable curves $\MM_{0,3}$ is also a point.
Thus
\[A^*(\mm_{0,0}(\PP^1,1))\iso\QQ \text{\hspace{.5in}and\hspace{.5in}}
A^*(\mm_{0,3}(\PP^1,0))\iso\QQ[H]/(H^2)\text{.}\]
In the latter, $H$ corresponds to the hyperplane pullback under any of the three
evaluation morphisms, which all simply 
record the image of the trivial stable map.

We also have 
$\mm_{0,1}(\PP^1,1)\iso \PP^1$
and $\mm_{0,2}(\PP^1,1)\iso \PP^1\x\PP^1$, where these moduli spaces simply parametrize the
images of the marked points. Proofs of these well-known facts are omitted here, but can
be found in \cite{C4}.
It follows that 
\[A^*(\mm_{0,1}(\PP^1,1))\iso\frac{\QQ[H_1]}{(H_1^2)}\text{.}\]
and
\[A^*(\mm_{0,2}(\PP^1,1))\iso\frac{\QQ[H_1,H_2]}{(H_1^2,H_2^2)}\text{.}\]

\begin{prop}
\label{m0311}
We have an isomorphism of $\QQ$-algebras 

\begin{eqnarray*}
& & A^*(\mm_{0,3}(\PP^1,1)) \\
& \iso & \frac{\QQ[H_1,H_2,H_3,D]}{(H_1^2,H_2^2,H_3^2, 
(H_1+H_2-D)(H_2+H_3-D), D(H_1-H_2), D(H_2-H_3))}
\end{eqnarray*}

\end{prop}

\noindent
{\bf Proof.}  Since $\mm_{0,3}(\PP^1,1)$ is the universal curve over
$\mm_{0,2}(\PP^1,1)$, it is easy to show that
$\mm_{0,3}(\PP^1,1) \iso\Bl_\D{(\PP^1)^3}$, where $\D$ is the small diagonal.
Now $\D\iso\PP^1$, and the restriction map 
corresponding to $i:\D\hookrightarrow(\PP^1)^3$ sends each $H_i$
to the hyperplane class in $\D$. So $i^*:A^* ((\PP^1)^3)\rightarrow
A^*(\D)$ is surjective, and we can take $H_1-H_2$ and $H_2-H_3$ as
generators for $\ker{i^*}$. The small diagonal is the complete
intersection of any two of the large diagonals. So we may apply Keel's
Lemma 1 from
\cite{K}, which says that whenever $X$ is a complete intersection of two
divisors $D_1$,$D_2$ in a scheme $Y$ and the restriction map 
$i^*:A^*(Y)\ra A^*(X)$ is
surjective, then
\[A^*(\tilde{Y})=A^*(Y)[T]/((D_1-T)(D_2-T), \ker{i^*}\cdot T)\text{,}\]
where $\tilde{Y}$ is the blowup of $Y$ along $X$. 
Here $T$ corresponds to the exceptional divisor. We know from \cite{F} that 
$A^*((\PP^1)^3)=\QQ[H_1,H_2,H_3]/(H_1^2,H_2^2,H_3^2)$ and that 
we can express two of the 
large 
diagonal classes as $D_i=H_i+H_{i+1}$ for $i\in\{1,2\}$.  
The expression in the proposition results.$\Box$

\noindent
The $H_i$ in the presentation are naturally identified with 
the corresponding hyperplane pullbacks. Furthermore $D$ corresponds to
the boundary divisor $D=D_{\und{3},0,\emptyset,1}$ whose generic stable map
has all three marked points
lying on the same collapsed component.

In \cite{MM}, Musta\c{t}\v{a} and Musta\c{t}\v{a} describe presentations 
for the Chow rings $A^*(\mm_{0,1}(\PP^r,2))$ as an illustration
of their more general results. Specifically, they find that

\[A^*(\mm_{0,1}(\PP^r,2))\iso\frac{\QQ[H,\psi,S,P]}{\I},\]
where the ideal $\I$ is generated by six relations. Setting $r=1$ and simplifying, 
we arrive at

\[I=(H^2, 3\psi^2H+\psi^3, P\psi, S(2H\psi+\psi^2), S(2H+3\psi)+S^2-2P, 4H+4\psi+S).\]

Here $H=H_1$ is the hyperplane divisor and $\psi=\psi_1$ is the $\psi$-class, just as described
in Section \ref{sec:gen}. Furthermore, in the framework of \cite{MM} $S$ corresponds
to the negative of the lone boundary class $D=D_{\underline{1},1,\emptyset,1}$. 
Hence $S=-D$ in our notation.
From the last relation, $\psi=-\frac{1}{4}S-H=\frac{1}{4}D-H_1$. This identity will be
important later. For now, we use it to further simplify $\I$ to

\[I=(H^2, S^3, -\frac{1}{4}P-PH, -HS+\frac{1}{4}S^2-2P)\]
where we have discarded $\psi$ as a generator. We can also discard $P$ and simplify
one more time using $P=\frac{1}{8}S^2-\frac{1}{2}HS$. We get

\[I=(H^2, S^3).\]

Thus, in our notation,

\[A^*(\mm_{0,1}(\PP^1,2))=\frac{\QQ[D,H_1]}{(H_1^2,D^3)}\text{.}\]

We noted above the crucial description $\psi=\frac{1}{4}D-H_1$ of the $\psi$-class in
$A^*(\mm_{0,1}(\PP^1,2))$ in terms of its boundary
and hyperplane divisors. We also need, and will now provide,
such expressions for $\psi$-classes on the moduli spaces $\mm_{0,n}(\PP^1,1)$, $n\in\und{3}$.
This is possible since we have seen that 
the the boundary and hyperplane divisor
classes generate the Chow rings in the cases under consideration. 
We will use the setup of 
Section \ref{sec:gen}.

To avoid confusion, we will use the notation $\pi_n$ for universal projections to
a space of stable maps (or more generally for contraction morphisms that 
forget a marked point) and $\r_i:(\PP^1)^n\rightarrow\PP^1$ for projections
of $(\PP^1)^n$ to a factor.

We will make use of several identities expressing pullbacks of
the standard divisor classes on $\m0$ under contraction morphisms in terms of
the standard divisor classes on $\mm_{0,n+1}(\PP^r,d)$.
The formula governing pullback of $\psi$-classes is 
$\psi_i=\pi_{n+1}^*(\psi_i)+D_{i,n+1}$, which is a well-known
extension of an identity in \cite{Wit}. On the
left side, $\psi_i\in A^*(\mm_{0,n+1}(\PP^r,d))$, while on the right side,
$\psi_i\in A^*(\m0)$ and $D_{i,n+1}$ is the divisor

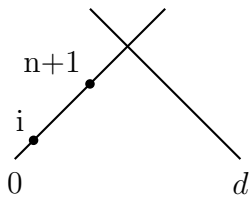
\begin{center}
\begin{pspicture}(0,0)(4,3)
\pnode(0.5,.5){a}
\dotnode(1.5,1.5){c}
\dotnode(.75,.75){d}
\pnode(2.5,2.5){e}
\pnode(1.5,2.5){f}
\pnode(3.5,0.5){j}
\ncline{a}{e}
\ncline{f}{j}
\uput{5pt}[d](.5,.5){0}
\uput{5pt}[d](3.5,.5){$d$}
\uput{5pt}[ul](.75,.75){i}
\uput{5pt}[ul](1.5,1.5){n+1}
\end{pspicture}
\end{center}
with the remaining marked points on the degree $d$ component. For $i<j$,
it's easy to see that $\ev_{i,n}\circ\pi_j=\ev_{i,n+1}$. Here the second
subscript of $\ev_i$ indicates the number of marked points associated to its
domain. It follows that $\pi_j^*(H_i)=H_i$ for $i<j$. A similar statement
holds when $i>j$ (and even when $i=j$ if the marked points are 
monotonically relabeled with $\und{n}$ in the target moduli space), 
but attention must be given to how the indexing changes in the
wake of deleting one index.

Finally, we recall the basic fact about pullbacks of boundary divisors: 
If $\pi$ is any contraction morphism, then

\begin{equation}
\label{pulldiv}
\pi^*(D_{A,d_A,B,d_B})=\sum_{A\sub A\pr,B\sub B\pr} D_{A\pr,d_A,B\pr,d_B}
\text{.}
\end{equation}
Clearly the righthand side is the support of the pullback.
See \cite{FP}, for example, for relevant statements about why the pullback is
multiplicity-free.

First we will consider the case $n=1$.
Recall that $\mm_{0,1}(\PP^1,1)\iso\PP^1$ and 
its universal curve is $\mm_{0,2}(\PP^1,1)\iso\PP^1\x\PP^1$.
Furthermore,
the section is the diagonal map $\D$. The universal 
projection is the projection $\r_1$.
Then 

\begin{eqnarray*}
\w_{\r_1} = K_{\PP^1\x\PP^1}\*\r_1^*K_{\PP^1}^\vee
=\r_1^*(\O(-2))\*\r_2^*(\O(-2))\*\r_1^*(\O(-2))^\vee=\O(0,-2)
\text{.}\end{eqnarray*}
Now $c_1(\O(0,-2))=-2H_2$, and the pullback of each $H_i$ under $\D$ is $H_1$.
So in this case $\psi=\psi_1=\D^*(c_1(\O(0,-2)))=-2H_1$.

Moving on to the case $n=2$ and
using the results stated above, in $A^*(\mm_{0,2}(\PP^1,1))$ we have
\[\psi_1=\r_1^*(-2H_1)+D_{1,2}=-2H_1+H_1+H_2=H_2-H_1\text{.}\]
Here $D_{1,2}$ is the divisor corresponding to the locus where the marked
points have the same image. This is the diagonal of $\PP^1\x\PP^1$, and we
have used the fact that its class in $A^*(\PP^1\x\PP^1)$ is $H_1+H_2$.

By symmetry, that is, by forgetting the first marked point via the
projection $\r_2$ instead of the second using $\r_1$,
we find
\[\psi_2=H_1-H_2\text{.}\]

Finally, for $n=3$, pulling back from $\mm_{0,2}(\PP^1,1)$ we have
\[\psi_1=\pi_3^*(\psi_1)+D_{1,3}=H_2-H_1+H_1+H_3-D=H_2+H_3-D
\text{,}\]
\[\psi_2=\pi_3^*(\psi_2)+D_{2,3}=H_1-H_2+H_2+H_3-D=H_1+H_3-D
\text{,}\]
and, by symmetry,
\[\psi_3=H_1+H_2-D
\text{.}\]
Here $D_{i,j}$ is the divisor corresponding to the closure of locus where the 
$i$'th and $j$'th marked
points have the same image, but the image of the remaining marked point
is different. This is the 
proper transform in $\Bl_\D{(\PP^1)^3}$ of the large
diagonal $\D_{ij}$ in $(\PP^1)^3$.
Its class in $A^*(\Bl_\D{(\PP^1)^3})$ is $H_i+H_j-D$.

\section{Relations}
\label{sec:rel}

Relations come from three different sources: Pullbacks of relations
in $A^*(\mm_{0,1}(\PP^r,2))$ via the contractions forgetting a marked
point, relations given by the geometry of $\psi$-classes
on boundary divisors, and one linear
relation that so far has no {\em proven} geometric explanation. Instead,
we prove this last relation using the method of localization and linear
algebra, which will be explained in Section \ref{sec:lla}.

\subsection{Relations from pullbacks}
\label{sec:pb}

In Section \ref{sec:simpler}, we found the relations $H_1^2$ and $D^3$ in
$A^*(\mm_{0,1}(\PP^1,2))$, as well as 
the relation $\psi-\frac{1}{4}D+H_1$.

There are two contraction morphisms $\pi_1$ and $\pi_2$
from $\mm_{0,2}(\PP^1,2)$ to
$\mm_{0,1}(\PP^1,2)$. Recall that $\pi_i$ forgets the $i$'th marked point. 
Formulas for the pullbacks of the standard divisor classes under these
maps were given in Section \ref{sec:simpler}.

First, pulling back the relation $H_1^2$ under $\pi_1$ and $\pi_2$
gives the relations $H_i^2$ for $i\in\und{2}$. (When considering the 
contraction $\pi_1$, either a monotonic relabeling of the marked points 
gives an index shift on pullback, or $H_1$ needs to be labeled as $H_2$.)
We can also see these relations via pullback under the evaluation morphisms.
Since $H^2=0$ in $A^*(\PP^1)$,
we have $H_i^2=\ev_i^*(H^2)=0$.

Second, applying Equation \ref{pulldiv} to the case at hand, 
we have $\pi_i^*(D)=D_1+D_2$ for $i\in\und{2}$.
Pulling back the relation $D^3$ in $A^*(\mm_{0,1}(\PP^1,2))$ by either one of 
these gives the cubic relation $(D_1+D_2)^3$ in $A^*(\mm_{0,2}(\PP^1,2))$.

Finally, 
linear relations expressing the $\psi$-classes in terms of the other divisor
classes are obtained by pulling back the relation $\psi-\frac{1}{4}D+H_1$.
We find
\[\psi_i=\pi_j^*(\psi_i)+D_0=\pi_j^*(\frac{1}{4}D-H_1)+D_0
=\frac{1}{4}D_1+\frac{1}{4}D_2+D_0-H_i\text{,}\]
where $i\neq j$. Thus we have relations
$\psi_i-\frac{1}{4}D_1-\frac{1}{4}D_2-D_0+H_i$.

\subsection{Relations from the geometry of $\psi$-classes on boundary 
divisors}
\label{sec:geomrel}

We can interpret a product with
a factor of $D_i$ as a
restriction to that divisor and, via gluing morphisms, ultimately as a class on
a fiber product of simpler 
moduli spaces.  This simplifies the computations of such products
once we know the pullbacks of the remaining factors.
We will use the gluing morphisms
\[j_0:\mm_{0,3}(\PP^1,0)\x_{\PP^1}\mm_{0,1}(\PP^1,2)
\rightarrow\mm_{0,2}(\PP^1,2)\]
and
\[j_1:\mm_{0,3}(\PP^1,1)\x_{\PP^1}\mm_{0,1}(\PP^1,1)
\rightarrow\mm_{0,2}(\PP^1,2)\text{,}\]
whose images are $D_0$ and $D_1$ respectively. By convention, 
both $j_0$ and $j_1$ glue the third
marked point of the first factor to the lone marked point of the second
factor. It is not hard to see that
$j_0$ is an isomorphism onto $D_0$, and we note further that
$\mm_{0,3}(\PP^1,0)\x_{\PP^1}\mm_{0,1}(\PP^1,2)\iso\mm_{0,1}(\PP^1,2)$.
Similarly, $\mm_{0,3}(\PP^1,1)\x_{\PP^1}\mm_{0,1}(\PP^1,1))
\iso\mm_{0,3}(\PP^1,1)$. However, $j_1$ is only an isomorphism away from
the divisor $D\x_{\PP^1}\mm_{0,1}(\PP^1,1)$, where $D$ is the boundary
divisor of $\mm_{0,3}(\PP^1,1)$ as described in Section \ref{sec:simpler}.
The image of this divisor is isomorphic to the global quotient stack 
\[[\mm_{0,1}(\PP^1,1)\x_{\PP^1}\mm_{0,4}(\PP^1,0)\x_{\PP^1}\mm_{0,1}(\PP^1,1)
/S_2]\text{,}\]
where the $S_2$-action switches the factors on the ends.
Thus the restriction of $j_0$ to $D$ has degree two.
Note that this last fiber product is isomorphic to 
$[\mm_{0,4}(\PP^1,0)/S_2]$,
where the $S_2$ action switches the third and fourth marked points.
 
For the $D_0$ case, the universal property of the moduli space shows that
$\psi_1$ and $\psi_2$ pull back to what may be considered as the first and
second $\psi$-classes on $\mm_{0,3}(\PP^1,0)$.
A similar statement holds for $D_1$ case,  
this time with the resulting $\psi$-classes on $\mm_{0,3}(\PP^1,1)$.

In the first case using this technique, we will show that the $\psi$-classes
vanish on $D_0$. This is because each 
$\psi$-class pulls back to zero under $j_0$ since the marked points lie
on the rigid component corresponding to $\mm_{0,3}(\PP^1,0)$. 
More rigorously, $j_0$ induces a family 
$(\C,\tilde{\s}_1,\tilde{\s}_2,\ev_3\circ\tilde{\jmath}_0)$ of
stable maps via the fiber diagram

\vspace{0.3in}

\psset{arrows=->}
\begin{center}
\begin{psmatrix}
$\C$  & $\mm_{0,3}(\PP^1,2)$ & $\PP^1$\\ 
$\mm_{0,3}(\PP^1,0)\x_{\PP^1}\mm_{0,1}(\PP^1,2)$ & 
$\mm_{0,2}(\PP^1,2)\text{,}$
\ncline{1,1}{2,1}\naput{$\tilde{\pi}$}
\ncline{1,1}{1,2}\naput{$\tilde{\jmath}_0$}
\ncline{1,2}{2,2}\naput{$\pi$}
\ncline{2,1}{2,2}\naput{$j_0$}
\ncline{1,2}{1,3}\naput{$\ev_3$}
\ncarc[arcangle=45]{2,1}{1,1}\naput{$\tilde{\s}_i$}
\ncarc[arcangle=45]{2,2}{1,2}\naput{$\s_i$}
\end{psmatrix}
\end{center}
\psset{arrows=-}

\vspace{0.2in}

\noindent
It is easy
to check that this family can be identified with a universal family of
stable maps over $D_0$. We can obtain another family over 
$\mm_{0,3}(\PP^1,0)\x_{\PP^1}\mm_{0,1}(\PP^1,2)$ by fiber product:

\vspace{0.3in}

\psset{arrows=->}
\begin{center}
\begin{psmatrix}
$\mm_{0,4}(\PP^1,0)\x_{\PP^1}\mm_{0,1}(\PP^1,2)$  & 
$\mm_{0,4}(\PP^1,0)$ & $\PP^1$ \\ 
$\mm_{0,3}(\PP^1,0)\x_{\PP^1}\mm_{0,1}(\PP^1,2)$ & 
$\mm_{0,3}(\PP^1,0)\text{,}$
\ncline{1,1}{2,1}\naput{$\tilde{\pi}_4$}
\ncline{1,1}{1,2}\naput{$\tilde{\pj}_1$}
\ncline{1,2}{2,2}\naput{$\pi_4$}
\ncline{2,1}{2,2}\naput{$\pj_1$}
\ncline{1,2}{1,3}\naput{$\ev_4$}
\ncarc[arcangle=45]{2,1}{1,1}\naput{$\tilde{s}_i$}
\ncarc[arcangle=45]{2,2}{1,2}\naput{$s_i$}
\end{psmatrix}
\end{center}
\psset{arrows=-}

\vspace{0.2in}

\noindent
where the right-hand part is the universal stable map.
Recall from Section \ref{sec:simpler} that the marked points in
$\mm_{0,3}(\PP^1,0)$ cannot vary.
It follows that $c_1(s_i^*(\w_{\pi_4}))=\psi_i=0$ for 
$i\in\und{3}$. Note also that $\w_{\pi_4}$
pulls back to the relative dualizing sheaf of the left-hand
column.
There is an inclusion 
\[k_0:\mm_{0,4}(\PP^1,0)\x_{\PP^1}\mm_{0,1}(\PP^1,2)\ra
\mm_{0,3}(\PP^1,2)\]
very similar to the morphism $j_0$ described above. This morphism is part
of a 2-commutative diagram that guarantees the existence of a stack morphism
\[\iota:\mm_{0,4}(\PP^1,0)\x_{\PP^1}\mm_{0,1}(\PP^1,2)\ra\C\]
over $\mm_{0,3}(\PP^1,0)\x_{\PP^1}\mm_{0,1}(\PP^1,2)$. Furthermore,
$\iota$ is injective since $k_0$ is, and 
it is also compatible with the sections.
The images
of the sections $\tilde{\s_i}$ are contained in the image of $\iota$.
Therefore we can compute the $\psi$-classes of
$\mm_{0,3}(\PP^1,0)\x_{\PP^1}\mm_{0,1}(\PP^1,2)$ on this subfamily:
\begin{eqnarray*}
\psi_i & = & c_1(\tilde{\s}_i^*\w_{\tilde{\pi}})\\
 & = & c_1(\tilde{s}_i^*\w_{\tilde{\pi}_4})\\
 & = & \tilde{s}_i^*\tilde{\pj}_1^*(c_1(\w_{\pi_4}))\\
 & = & \pj_1^*s_i^*(c_1(\w_{\pi_4}))\\
 & = & 0\text{.}
\end{eqnarray*}
It follows that the pullbacks of $\psi_1$ and $\psi_2$
to $D_0$ vanish. This gives relations $D_0\psi_i$.

Now we will show that the product $D_1\psi_1\psi_2$ vanishes by computing the
pullback of $\psi_1\psi_2$ under $j_1$. 
Identifying $A^*(\mm_{0,3}(\PP^1,1)\x_{\PP^1}\mm_{0,1}(\PP^1,1))$ with
$A^*(\mm_{0,3}(\PP^1,1))$, the pullback is
\[\psi_1\psi_2=(H_2+H_3-D)(H_1+H_3-D)
=(H_1+H_2-D)(H_2+H_3-D)=0\]
according to the presentation given in Proposition \ref{m0311}. 
This shows that $\psi_1\psi_2$ is in the kernel of $j_1^*$.
It follows from \cite[Theorem 3.1]{Mu} that there is an group isomorphism
$A^2(D_1)\ra A^2(\mm_{0,3}(\PP^1,1))^{S_2}$, where the target is the subgroup
of $S_2$-invariants. This isomorphism is naturally identified with a
monomorphism into $A^2(\mm_{0,3}(\PP^1,1))$.
Furthermore, Mumford's proof of this theorem 
shows that, up to sign, this monomorphism
is the same as $\tilde{\jmath}_1^*$, where $\tilde{\jmath}_1$ is 
$j_1$ with the target
changed to $D_1$ (and considered as a map of these 
Chow groups). This is enough
to show that $\psi_1\psi_2$ restricts to zero on $D_1$.
Thus $D_1\psi_1\psi_2$ is a relation in $A^*(\mm_{0,2}(\PP^1,2))$.

\subsection{A relation from localization and linear algebra}
\label{sec:lla}

Ideally, the proof that (\ref{geomprez}) is a presentation for 
$A^*(\mm_{0,2}(\PP^1,2))$ should be uniformly geometric. Such a proof
would be easier to understand and generalize to higher-dimensional targets
$\PP^r$ (and other moduli spaces). 
Nonetheless, we will use a sort of ``brute force" algebraic computation
for one aspect of the proof: the
derivation of the linear relation $D_2-\psi_1-\psi_2$. It is evident that
the relation has a geometric source. Indeed, the heuristic calculations of 
\cite{Wit} used to justify relations on moduli spaces of stable curves
lead us to believe that, similarly, $D_2-\psi_1-\psi_2$ should be explained 
geometrically by the existence of a certain section of the tensor product
$L_1\* L_2$ of the cotangent line bundles. 

\begin{conj}
There exists a section $s$ of the line bundle $L_1\* L_2$ on $A^*(\mm_{0,2}(\PP^1,2))$
whose zero stack is the boundary divisor $D_2$.
\end{conj}
A rigorous construction of such $s$ would give the relation, for then
\[D_2=Z(s)=c_1(L_1\* L_2)=c_1(L_1)+c_1(L_2)=\psi_1+\psi_2\text{.}\]
%
%Lacking such a construction, we proceed to a less elegant computational
%approach, which we hope is of interest in its own way.

\begin{rmk}
D. Oprea has noted that, alternatively, this relation follows by adding the relations
$H_1=H_2+2\psi_2-D_2$ and $H_2=H_1+2\psi_1-D_2$ first written down by
Lee and Pandharipande in \cite{LP}. He further communicated the 
following geometric explanation of these relations. The $\psi$-class
$\psi_1$ on $\mm_{0,3}(\PP^1,2)$ can be expressed as both 
$\psi_1=D_{\{1\},1,\{2,3\},1}+D_{\{1\},2,\{2,3\},0}$ and
$\psi_1=\pi_3^*(\psi_1)+D_{\{2\},2,\{1,3\},0}$, using the notation of
Section \ref{sec:gen}. Equating these, intersecting
with $H_3$, and pushing forward by $\pi_3$ gives the second relation
above. The computation for the first relation is symmetric.
\end{rmk}

We call the computational algebraic method developed in this subsection
{\em localization and linear algebra}. In theory, it could be used to compute
relations in any moduli space of stable maps to projective space. (In 
practice, the method quickly becomes too tedious as the parameters increase. See 
\cite[Section 7.3]{C4} for an illustration.)
Its first step consists of using
localization to
find the integrals of all degree four monomials in the boundary divisors and
hyperplane pullbacks in $A^*(\mm_{0,2}(\PP^1,2))$.
Then an arbitrary relation of given degree in these generating classes is
considered (with variable coefficients). 
This relation is multiplied by various monomials of complementary
dimension, and the result is integrated. The resulting rational polynomials 
in the coefficients place restrictions on these coefficients. Ultimately,
we can use these restrictions to describe precisely the form of 
relations that can occur in each degree in this Chow ring.

Before proceeding, we briefly recall the requisite localization results.
Localization was first applied to the spaces $\m0$ by Kontsevich in \cite{Ko}.
We will follow the exposition of \cite[Chapter 9]{CK}.

Let $T=(\CC^*)^n$, and let $H_T^*(X)$ denote the $T$-equivariant cohomology ring of $X$.
(This is isomorphic to the $T$-equivariant Chow ring of \cite{EG} in our case.)
The ring $H_T^*(X)$ includes generators $\la_1, \ldots ,\la_n$, where $\la_i$ is the weight
of the character given by the $i$'th projection $(\CC^*)^n\ra\CC^*$.  
Let $\R_T\iso\CC(\la_1,\ldots,\la_n)$. The fixed point locus $X^T$ is a
union of smooth connected components $Z_j$. We have inclusions
$i_j:Z_j\ra X$ and normal bundles $N_j=N_{Z_j/X}$ which are equivariant.
Finally, let $\eqeul(N_j)$ be the equivariant Euler class of the normal
bundle for each $j$. The residue formula that makes our computations possible
is the following corollary of the localization theorem of Atiyah-Bott.

\begin{thm}
\label{stackloc}
Let $X$ be an orbifold which is the variety underlying a smooth
stack with a $T$-action. 
If $\a\in H_T^*(X)\*\R_T$, then
\[\int_{X_T}\a=\sum_j \int_{(Z_j)_T}\frac{i_j^*(\a)}{a_j\eqeul(N_j)}\text{,}\]
where $a_j$ is the order of the group $H$ occurring in a local chart at the
generic point of $Z_j$.
\end{thm}

The natural action of 
$T=(\CC^*)^{r+1}$ on $\PP^r$ induces a $T$-action
on the smooth Deligne-Mumford stack $\m0$ by composition of the action with stable maps. 
Let $q_0,\ldots, q_r$ be the fixed points of $\PP^r$ 
under this $T$-action.
A $T$-fixed point of $\m0$ corresponds
to a stable map $(C, p_1,\ldots, p_n, f)$ where each component of $C_i$ of $C$
is either mapped to some $q_i$ or multiply covers a 
coordinate line.  Each marked point $p_j$, each node of $C$, and each
ramification point of $f$ is mapped to a $q_i$ also. 
This gives a 1--1 
correspondence of the connected components of $\m0^T$ with connected trees $\G$ of the following type:
The vertices $v$ of $\G$ are in 1--1 correspondence with the connected 
components $C_v$ of $f^{-1}(\{q_0,\ldots,q_r\})$, so each $C_v$ is either a point
or a connected union of irreducible components of $C$. The edges
$e$ of $\G$ correspond to irreducible components $C_e$ of $C$ which are
mapped onto some coordinate line $\ell_e$ in $\PP^r$.

The graph $\G$ has the following labels: Associate to each vertex $v$
the number $i_v$ defined by $f(C_v)=q_{i_v}$, as well as the set $S_v$
consisting of those $i$ for which the marked point $p_i$ is in $C_v$.
Associate to each edge $e$ the degree $d_e$ of the map $f|_{C_e}$.
Finally, we impose the following three conditions:

\begin{enumerate}
\item If an edge $e$ connects $v$ and $v\pr$, then $i_v\neq i_{v\pr}$, and
$\ell_e$ is the coordinate line joining $q_{i_v}$ and $q_{i_{v\pr}}$.

\item $\sum_e d_e=d$.

\item $\coprod_v S_v=\und{n}$.
\end{enumerate}
 
Define $\MM_\G$ to be $\prod_{v:\dim C_v=1} \MM_{0,n(v)}$. It has a
group of automorphisms $A_\G$ whose order is $a_{\G}=|\Aut(\G)|\prod_e{d_e}$,
where $\Aut(\G)$ is the group of automorphisms of $\G$ which preserve
the labels. The locus of stable maps with graph $\G$ is the quotient stack $[\MM_\G/A_\G]$. 
When Theorem \ref{stackloc} is applied to $\m0$, the scalar factor
appearing in the denominator of the term corresponding to this fixed
component is $a_{\G}$.

The last ingredients needed in order to use localization on $\m0$ are the
Euler classes of the fixed components. Denote the normal bundle of
$\mm_\G$ by $N_\G$. Define a {\em flag} $F$ of a graph to be a pair
$(v,e)$ such that $v$ is a vertex of $e$. Put $i(F)=v$ and let $j(F)$ be
the other vertex of $e$. Set
\[\w_F=\frac{\la_{i_{i(F)}}-\la_{i_{j(F)}}}{d_e}\text{.}\]
This corresponds to the weight of the $T$-action on the tangent space
of the component $C_e$ of $C$ at the point $p_F$ lying over $i_v$. 
%(In \cite[Chapter 9]{CK}, the Euler class of the
%tangent space to $\PP^r$ at $q_i$ is shown to be 
%
%\[\prod_{j\neq i} (\la_i-\la_j)\text{,}\]
%
%and the degree $d_e$ of the map introduces a denominator upon pullback.)
Let
$e_F$ be the first Chern class of the bundle on $\mm_\G$ whose fiber
is the cotangent space to the component associated to $v$ at $p_F$.
(More information about this type of class was given in Section \ref{sec:gen}.)
If $\val(v)=1$, let $F(v)$ denote
the unique flag containing $v$.  If $\val(v)=2$, let $F_1(v)$ and $F_2(v)$
denote the two flags containing $v$. Similarly, let $v_1(e)$ and $v_2(e)$
be the two vertices of an edge $e$.

\begin{thm}\label{norm}
The equivariant Euler class of the normal bundle $N_\G$ is a product
of contributions from the flags, vertices and edges:
\[\eqeul(N_\G)=e_{\G}^{\text{F}}e_{\G}^{\text{v}}e_{\G}^{\text{e}}\text{,}\]
where 
\[e_{\G}^{\text{F}}=\frac{\prod_{F:n(i(F))\geq 3}(\w_F-e_F)}
{\prod_{F}\prod_{j\neq i_{i(F)}}(\la_{i_{i(F)}}-\la_j)}\]
\[e_{\G}^{\text{v}}=\left(\prod_v \prod_{j\neq i_v}(\la_{i_v}-\la_j)
\right)
\left(\prod_{{\val(v)=2\atop S_v=\emptyset}}(\w_{F_1(v)}+\w_{F_2(v)})\right)
/\prod_{{\val(v)=1\atop S_v=\emptyset}} \w_{F(v)}\]
\[e_{\G}^{\text{e}}=\prod_e\left(\frac{(-1)^{d_e}(d_e!)^2
(\la_{i_{v_1(e)}}-\la_{i_{v_2(e)}})^{2d_e}}{d_e^{2d_e}}
\prod_{{a+b=d_e\atop k\neq i_{v_j(e)}}}
\left(\frac{a\la_{i_{v_1(e)}}+b\la_{i_{v_2(e)}}}{d_e}-\la_k\right)\right)
\text{.}
\]
\end{thm}

We now have the machinery necessary to compute the integrals of degree four monomials in
$A^*(\mm_{0,2}(\PP^1,2))$. Notice first
that
the relations $D_0\psi_i$ and $\psi_i-\frac{1}{4}D_1-\frac{1}{4}D_2-D_0+H_i$
imply 
\begin{equation}\label{d0hi}
D_0H_1=D_0H_2\text{,}
\end{equation} 
and hence $H_1H_2D_0=H_1^2D_0=0$.
The integrals of the following types of monomials will be zero since the 
monomials themselves are zero:

\begin{enumerate}

\item Any monomial with a factor of $H_i^2$ for $i\in\underline{2}$.

\item Any monomial with a factor of $H_1H_2D_0$.

\end{enumerate}

These together take care of $2\left(\left({5 \atop 1}\right)+
\left({5 \atop 2}\right)\right)-1+3=32$ of the $\left({8\atop 4}\right)=70$ 
degree 4 monomials. We will use localization to compute the remaining 38 
integrals.  (While not all of these values are necessary to obtain our relation,
we will list them all for completeness.) Let $T=(\CC^*)^2$. 
Consider the usual $T$-action on $\mm$.  First we have to find the 
$T$-equivariant 
Euler classes of the normal bundles of the fixed point components.
We will label the graphs 
of the 14 fixed components as follows.

\[
\begin{array}{cc}
\G_1=\text{
\psset{labelsep=2pt, tnpos=b,radius=2pt}
\pstree[treemode=R]{\TC*~{0}~[tnpos=a]{\{1\}}}
{
\TC*~{1}~[tnpos=a]{\{2\}}\taput{2} 
}} &
\G_8=\text{
\psset{labelsep=2pt, tnpos=b,radius=2pt}
\pstree[treemode=R]{\TC*~{0}~[tnpos=a]{\{1\}}}
{
\pstree{\TC*~{1}~[tnpos=a]{\{2\}}\taput{1}}
{
\TC*~{0}\taput{1} 
}
}
}\vspace{0.2in}
\\ \vspace{0.2in}

\G_2=\text{
\psset{labelsep=2pt, tnpos=b,radius=2pt}
\pstree[treemode=R]{\TC*~{0}~[tnpos=a]{\{2\}}}
{
\TC*~{1}~[tnpos=a]{\{1\}}\taput{2} 
}
} &
\G_9=\text{
\psset{labelsep=2pt, tnpos=b,radius=2pt}
\pstree[treemode=R]{\TC*~{1}~[tnpos=a]{\{2\}}}
{
\pstree{\TC*~{0}~[tnpos=a]{\{1\}}\taput{1}}
{
\TC*~{1}\taput{1} 
}
}
}

\\ \vspace{0.2in}

\G_3=\text{
\psset{labelsep=2pt, tnpos=b,radius=2pt}
\pstree[treemode=R]{\TC*~{0}~[tnpos=a]{\underline{2}}}
{
\TC*~{1}\taput{2} 
}
}  &
\G_{10}=\text{
\psset{labelsep=2pt, tnpos=b,radius=2pt}
\pstree[treemode=R]{\TC*~{0}~[tnpos=a]{\{2\}}}
{
\pstree{\TC*~{1}~[tnpos=a]{\{1\}}\taput{1}}
{
\TC*~{0}\taput{1} 
}
}
}\\ \vspace{0.2in}

\G_4=\text{
\psset{labelsep=2pt, tnpos=b,radius=2pt}
\pstree[treemode=R]{\TC*~{0}}
{
\TC*~{1}~[tnpos=a]{\underline{2}}\taput{2} 
}
}  &
\G_{11}=\text{
\psset{labelsep=2pt, tnpos=b,radius=2pt}
\pstree[treemode=R]{\TC*~{1}}
{
\pstree{\TC*~{0}~[tnpos=a]{$\underline{2}$}\taput{1}}
{
\TC*~{1}\taput{1} 
}
}
}\\ \vspace{0.2in}

\G_5=\text{
\psset{labelsep=2pt, tnpos=b,radius=2pt}
\pstree[treemode=R]{\TC*~{1}~[tnpos=a]{$\underline{2}$}}
{
\pstree{\TC*~{0}\taput{1}}
{
\TC*~{1}\taput{1} 
}
}
} &
\G_{12}=\text{
\psset{labelsep=2pt, tnpos=b,radius=2pt}
\pstree[treemode=R]{\TC*~{0}}
{
\pstree{\TC*~{1}~[tnpos=a]{$\underline{2}$}\taput{1}}
{
\TC*~{0}\taput{1} 
}
}
}\\ \vspace{0.2in}

\G_6=\text{
\psset{labelsep=2pt, tnpos=b,radius=2pt}
\pstree[treemode=R]{\TC*~{0}~[tnpos=a]{$\underline{2}$}}
{
\pstree{\TC*~{1}\taput{1}}
{
\TC*~{0}\taput{1} 
}
}
} &
\G_{13}=\text{
\psset{labelsep=2pt, tnpos=b,radius=2pt}
\pstree[treemode=R]{\TC*~{1}~[tnpos=a]{\{1\}}}
{
\pstree{\TC*~{0}\taput{1}}
{
\TC*~{1}~[tnpos=a]{\{2\}}\taput{1} 
}
}
}\\ \vspace{0.2in}

\G_7=\text{
\psset{labelsep=2pt, tnpos=b,radius=2pt}
\pstree[treemode=R]{\TC*~{1}~[tnpos=a]{\{1\}}}
{
\pstree{\TC*~{0}~[tnpos=a]{\{2\}}\taput{1}}
{
\TC*~{1}\taput{1} 
}
}
} &
\G_{14}=\text{
\psset{labelsep=2pt, tnpos=b,radius=2pt}
\pstree[treemode=R]{\TC*~{0}~[tnpos=a]{\{1\}}}
{
\pstree{\TC*~{1}\taput{1}}
{
\TC*~{0}~[tnpos=a]{\{2\}}\taput{1} 
}
}
}
\end{array}
\]

\vspace{0.1in}

\noindent
Let $Z_i$ denote the fixed component corresponding to $\G_i$ for all $i$. 
All of the fixed components are points except for $Z_{11}$ and $Z_{12}$,
which are each isomorphic to the quotient $[\PP^1/S_2]$. Fixed components
$Z_1$, $Z_2$, $Z_3$, and $Z_4$ also have automorphism group $S_2$.

All
nonzero degree 4 classes have factors $D_i$, and thus are supported on the boundary of
$\mm_{0,2}(\PP^1,2)$. Hence representatives of these classes 
will not intersect $Z_1$ or $Z_2$, which lie in the locus 
where the domain curves are smooth.
Since restrictions of the relevant classes to $Z_1$ and $Z_2$ are thus zero,
their equivariant Euler classes are not needed for our computation. 

The term $e_F$ that appears in the formula for $e_{\G}^{\text{F}}$ in
Theorem \ref{norm} 
is clearly zero on any fixed component which is a point.  
It follows from the discussion of $\tilde{\psi}$-classes in
Section \ref{sec:gen} that,
for the fixed components that are isomorphic to $[\PP^1/S_2]$,
we can identify $e_F$ with the class of a point (mod 2).
See \cite{C4} for a proof.

Straightforward application of Theorem \ref{norm}, together with other facts stated
in \cite[Chapter 9]{CK}, leads to the following equivariant
Euler classes.

\[\eqeul(N_{\G_{3}})=\eqeul(N_{\G_{4}})
=\frac{-(\la_1-\la_0)^4}{4}
\]

\[\eqeul(N_{\G_{5}})=\eqeul(N_{\G_{6}})=2(\la_1-\la_0)^4
\]

\[\eqeul(N_{\G_{7}})=\eqeul(N_{\G_{8}})=\eqeul(N_{\G_{9}})=
\eqeul(N_{\G_{10}})=-(\la_0-\la_1)^4
\]

\[\eqeul(N_{\G_{11}})
=(\la_0-\la_1)^2(\la_0-\la_1-2\psi)
\]

\[\eqeul(N_{\G_{12}})=(\la_1-\la_0)^2(\la_1-\la_0-2\psi)
\]

\[\eqeul(N_{\G_{13}})=\eqeul(N_{\G_{14}})=2(\la_0-\la_1)^4
\]

Next we need to know the restriction of each degree 4 monomial in the 
generating classes to each fixed component.  These come immediately from the
restrictions of the generating classes themselves to the fixed components,
which are given in the Table \ref{rest0212}. These values again follow from
the theory laid out in \cite[Chapter 9]{CK}. Some care must be taken when
deriving the nonzero restrictions, especially when there are automorphisms involved.
(Technically, Table \ref{rest0212} gives the pullbacks of the divisor classes to an
atlas of each fixed component; the effect of quotienting by the automorphism groups
is taken into account later during application of Theorem \ref{stackloc}.)
Detailed arguments are given in \cite{C4}.

%\vspace{0.1in}

\begin{table}
%\begin{center}
\begin{tabular*}{5.9375in}{|c|c|c|c|c|c|c|c|c|@{\extracolsep{\fill}}c|}
\hline
%{\small F. comp.} 
& $Z_1$ & $Z_2$ & $Z_3$ & $Z_4$ & $Z_5$ & $Z_6$ & $Z_7$
& $Z_8$  & $Z_9$\\ \hline 

$H_1$ & $\la_0$ & $\la_1$ &  $\la_0$ &  $\la_1$ &  $\la_1$ &  $\la_0$ &  $\la_1$
& $\la_0$ & $\la_0$  \\ \hline

$H_2$ &  $\la_1$ & $\la_0$ & $\la_0$ &  $\la_1$ & $\la_1$ & $\la_0$ & $\la_0$ 
 & $\la_1$ & $\la_1$ \\ \hline

$D_0$ & 0 & 0 & $\frac{\la_0-\la_1}{2}$ & $\frac{\la_1-\la_0}{2}$ 
& $\la_1-\la_0$ & $\la_0-\la_1$ & 0 & 0  & 0 \\ \hline

$D_1$ & 0 & 0 & 0 & 0 & $2(\la_0-\la_1)$ & $2(\la_1-\la_0)$ & 
$\la_0-\la_1$ & $\la_1-\la_0$  & $\la_0-\la_1$ \\ \hline

$D_2$ & 0 & 0 & 0 & 0 & 0 & 0 & $\la_0-\la_1$ & $\la_1-\la_0$ & $\la_0-\la_1$ 
\\ \hline

\end{tabular*}

\begin{tabular}{|c|c|c|c|c|c|c|}
\hline
%{\small F. comp.}  
& $Z_{10}$ & $Z_{11}$ & $Z_{12}$ & 
$Z_{13}$ & $Z_{14}$
 \\ \hline

$H_1$ &  $\la_1$ &  $\la_0$ &  $\la_1$ &  $\la_1$ &  $\la_0$
  \\ \hline

$H_2$  & $\la_0$ & $\la_0$ & $\la_1$ & $\la_1$ & $\la_0$
  \\ \hline

$D_0$ & 0 & $\psi$ & $\psi$ & 0 & 0   \\ \hline

$D_1$ & $\la_1-\la_0$ &  $2\la_0-2\la_1-2\psi$ 
& $2\la_1-2\la_0-2\psi$ & 0 & 0  \\ \hline

$D_2$  & $\la_1-\la_0$ &  $2\psi$ &  $2\psi$ & $2(\la_0-\la_1)$ & 
$2(\la_1-\la_0)$   \\ \hline

\end{tabular}

\caption{Restrictions of divisor classes in $A_T^*(\mm_{0,2}(\PP^1,2))$
to fixed components
\label{rest0212}}
%\end{center}
\end{table}
%\vspace{0.2in}

We now have everything we need to compute the integrals.  
Computations of all the integrals can be found in \cite[Appendix 1]{C4}, giving the results 
shown in Table \ref{deg4ints}. Any integral of a degree four monomial not listed there is automatically
zero for one of the reasons given at the beginning of the section.
Details are shown below for the integral of $D_1D_2H_1H_2$ to give a flavor of the calculations. 
Note first that
\[(\la_0-\la_1-2\psi)^{-1}=(\la_0-\la_1)^{-1}(1-2\psi/(\la_0-\la_1))^{-1}
=\frac{1+2\psi/(\la_0-\la_1)}{(\la_0-\la_1)}\text{,}\]
(since $\psi^2=0$) and similarly with $\la_0$ and $\la_1$ switched. Also note the factors of
two appearing in the denominators of the integrands for $Z_{11}$ and $Z_{12}$,
which are due to the $S_2$ automorphism group of these components. 
We have

\begin{eqnarray*}
\int_{\mm_T}D_1D_2H_1H_2
& = & \int_{(Z_{7})_T}\frac{\la_0\la_1(\la_0-\la_1)^2}{-(\la_0-\la_1)^4}
+  \int_{(Z_{8})_T}\frac{\la_0\la_1(\la_1-\la_0)^2}{-(\la_0-\la_1)^4} \\
&   & +\int_{(Z_{9})_T}\frac{\la_0\la_1(\la_0-\la_1)^2}{-(\la_0-\la_1)^4} 
+\int_{(Z_{10})_T}\frac{\la_0\la_1(\la_1-\la_0)^2}{-(\la_0-\la_1)^4} \\
&   & +\int_{(Z_{11})_T}\frac{\la_0^22\psi(2\la_0-2\la_1-2\psi)} 
{2(\la_0-\la_1)^2(\la_0-\la_1-2\psi)} \\
&   &
+\int_{(Z_{12})_T}\frac{\la_1^22\psi(2\la_1-2\la_0-2\psi)}
{2(\la_1-\la_0)^2(\la_1-\la_0-2\psi)} \\
& = & -4\frac{\la_0\la_1}{(\la_0-\la_1)^2}+
\int_{(Z_{11})_T}\frac{2\la_0^2\psi(2\la_0-2\la_1)(1+2\psi/(\la_0-\la_1))}
{2(\la_0-\la_1)^3} \\
&  & +\int_{(Z_{12})_T}\frac{2\la_1^2\psi(2\la_1-2\la_0)(1+2\psi/(\la_1-\la_0))}
{2(\la_1-\la_0)^3}\\
& = & \frac{-4\la_0\la_1}{(\la_0-\la_1)^2}+\int_{(Z_{11})_T}\frac{2\la_0^2\psi}
{(\la_0-\la_1)^2}+\int_{(Z_{12})_T}\frac{2\la_1^2\psi}{(\la_1-\la_0)^2} \\
& = & \frac{2\la_0^2-4\la_0\la_1+2\la_1^2}{(\la_0-\la_1)^2} \\
& = & 2
\text{.}\end{eqnarray*}

\begin{table}
\begin{center}
\begin{tabular}{|p{2in}p{2.8in}|}
\hline
\rule[-3mm]{0mm}{8mm}$\int_{\mm}D_2^4=12$ & $\int_{\mm}D_2^3H_1=-4$ \\ 

\rule[-3mm]{0mm}{8mm}$\int_{\mm}D_2^3D_1=-4$ & $\int_{\mm}D_2^3H_2=-4$ \\
\rule[-3mm]{0mm}{8mm}
$\int_{\mm}D_2^3D_0=0$ & $\int_{\mm}D_2^2D_1H_1=0$\\
\rule[-3mm]{0mm}{8mm}
$\int_{\mm}D_2^2D_1^2=-4$ & $\int_{\mm}D_2^2D_1H_2=0$ \\
\rule[-3mm]{0mm}{8mm}
$\int_{\mm}D_2^2D_1D_0=0$ & $\int_{\mm}D_2^2D_0H_1=\int_{\mm}D_2^2D_0H_2=0$ \\
\rule[-3mm]{0mm}{8mm}
$\int_{\mm}D_2^2D_0^2=0$ & $\int_{\mm}D_2D_1^2H_1=4$\\
\rule[-3mm]{0mm}{8mm}
$\int_{\mm}D_2D_1^3=12$  & $\int_{\mm}D_2D_1^2H_2=4$\\
\rule[-3mm]{0mm}{8mm}
$\int_{\mm}D_2D_1^2D_0=0$ & $\int_{\mm}D_2D_1D_0H_1=\int_{\mm}D_2D_1D_0H_2=0$ 
\\
\rule[-3mm]{0mm}{8mm}
$\int_{\mm}D_2D_1D_0^2=0$ &  $\int_{\mm}D_2D_0^2H_1=\int_{\mm}D_2D_0^2H_2=0$\\
\rule[-3mm]{0mm}{8mm}
$\int_{\mm}D_2D_0^3=0$ &  $\int_{\mm}D_1^3H_1=-8$\\
\rule[-3mm]{0mm}{8mm}
$\int_{\mm}D_1^4=-20$ & $\int_{\mm}D_1^3H_2=-8$ \\
\rule[-3mm]{0mm}{8mm}
$\int_{\mm}D_1^3D_0=0$  & $\int_{\mm}D_1^2D_0H_1=\int_{\mm}D_1^2D_0H_2=4$\\
\rule[-3mm]{0mm}{8mm}
$\int_{\mm}D_1^2D_0^2=4$ & $\int_{\mm}D_1D_0^2H_1=\int_{\mm}D_1D_0^2H_2=-1$ \\
\rule[-3mm]{0mm}{8mm}
$\int_{\mm}D_1D_0^3=-2$ & $\int_{\mm}D_0^3H_1=\int_{\mm}D_0^3H_2=\frac{1}{4}$
\\
\rule[-3mm]{0mm}{8mm}
$\int_{\mm}D_0^4=\frac{3}{4}$ & $\int_{\mm}D_2D_1H_1H_2=2$ \\
\rule[-3mm]{0mm}{8mm}
$\int_{\mm}D_2^2H_1H_2=2$ & $\int_{\mm}D_1^2H_1H_2=2$ \\ \hline
\end{tabular}
\end{center}
\caption{Integrals of degree four classes on $\mm_{0,2}(\PP^1,2)$
\label{deg4ints}}
\end{table}

\begin{table}
%\hspace{-1.7in}
\begin{center}
\begin{tabular}{|c|c|}
\hline
Monomial & Resulting relation on coefficients \\ \hline
 $D_2^2H_1$ & $ 2b-4e=0$ \\ 
 $D_2^2H_2$ & $ 2a-4e=0$ \\ 
 $D_1D_0H_1$ & $ -c+4d=0$ \\ 
$D_1H_1H_2$ & $ 2d+2e=0$ \\ \hline
\end{tabular}
\caption{Restrictions placed on coefficients of a linear relation by 
integration \label{rest}}
\end{center}
\end{table}

Now we will use these results to find the desired linear relation in $A^*(\mm_{0,2}(\PP^1,2))$.
Expressions for the $\psi$-classes in $\mm_{0,2}(\PP^1,2)$ were given
in Section \ref{sec:pb}, so it suffices to consider the
basic divisor classes
$H_1$, $H_2$, $D_0$, $D_1$, and $D_2$.
Since the first Betti number is four, there must be a relation among these
five divisor classes.  Suppose we have a relation
\[aH_1+bH_2+cD_0+dD_1+eD_2=0\text{.}\]
We can place restrictions on the coefficients by multiplying the above 
equation by degree three monomials and then integrating. 
For example, multiplying by $D_2^2H_1$ gives
\[aD_2^2H_1^2+bD_2^2H_2H_1+cD_2^2D_0H_1+dD_2^2D_1H_1+eD_2^3H_1=0\text{.}\]
Now integration gives the equation
$2b-4e=0$, using the integral values from Table \ref{deg4ints}. 
Continuing with some other choices of monomial, we 
get the system of 
restrictions given in Table \ref{rest}.

%\vspace{0.2in}

%\vspace{0.2in}

Together these restrictions show that up to a constant multiple
the only possible linear
relation among these five divisor classes is
\[2H_1+2H_2-4D_0-D_1+D_2=0\text{.}\]
Thus this must indeed be a relation, and the 
remaining four classes must be independent.  Hence, we have additionally
found that 
the classes $H_1$, $H_2$,
$D_0$, and $D_1$ generate the degree one piece of the graded ring
$A^*(\mm_{0,2}(\PP^1,2))$. 
%The same method will be used in Section
%\ref{sec:prez} to show monomials in the basic 
%classes generate in degrees two and three also.
We can write the relation above as 
$D_2=4D_0+D_1-2H_1-2H_2$. Notice that this linear relation 
can also be written $D_2-\psi_1-\psi_2=0$ using the expressions for the 
$\psi$-classes from Section \ref{sec:pb}.

An interesting consequence of this relation, together with the relations
$D_0\psi_i$, is the relation $D_0D_2=0$. We can find this relation directly
by arguing that the divisors $D_0$ and $D_2$ are disjoint.
We have seen that
the domain curve of every stable map lying in
$D_0$ must have both marked points on the same degree zero component, and that
the
domain curve of every
stable map in $D_2$ must have the marked points on distinct
components.
These mutually exclusive properties of stable maps in $D_0$ and $D_2$
validate the claim of their disjointness.  Indeed, one may need to use such
a direct argument to get this kind of relation in other cases,
although the relation $D_2-\psi_1-\psi_2$ holds in $A^*(\mm_{0,2}(\PP^r,2))$ 
for general $r$ by the results of Lee and Pandharipande mentioned earlier (\cite{LP}).

\section{The Presentation}
\label{sec:prez}

\subsection{
Completeness of the higher degree parts of the presentation}
\label{sec:comp}

The only question that remains is whether there could be further relations
among the generators we have given, offset by additional generators we have missed.
In this subsection, we confirm that there are no more independent generators or
relations beyond those we have identified.
We have already seen in both Sections \ref{sec:gen} and \ref{sec:lla} that there can be no extra 
divisor classes independent of the ones given in Section \ref{sec:gen}, since
the first Betti number is four, and four of these classes are independent.
The same method used in Section \ref{sec:lla} also demonstrates that monomials in the
standard divisor classes generate in higher degrees. (This could
be argued using the additive basis as well, but we will not do so here.)
As in Section \ref{sec:lla}, we need not consider the $\psi$-classes since
they can be expressed in terms of the other generators. 
We need not consider monomials involving $D_2$ for the same reason:
As found in Section \ref{sec:lla}, $D_2$ can be expressed in terms of
other generators. Furthermore,
relations involving the $\psi$-classes can be reformulated to reduce
the number of spanning monomials in the other divisor classes.

In degree 2, there are ten monomials in the remaining 
classes. However, the $H_i^2$
will not play a role since they vanish. Furthermore,
it is easy to check that $\psi_1-\psi_2=H_1-H_2$, so that
$D_0H_1-D_0H_2=D_0\psi_1-D_0\psi_2=0$. Hence we can also discount the
monomial $D_0H_2$, leaving seven monomials spanning the degree two part
of the Chow ring. Suppose we have a 
relation of the form 
\[aD_1^2+bD_1D_0+cD_1H_1+dD_1H_2+eD_0^2+fD_0H_1+gH_1H_2=0\text{.}\]
Then, multiplying this expression by each of these seven degree two monomials
and integrating the results, we obtain a system of seven linear equations
in seven variables.  Solving this system, we see that the only possible 
relation of this form (up to a constant multiple) is 
\[D_1D_0+4D_0^2-4D_0H_1=0\text{.}\]
Using the expression for $D_2$ derived in Section \ref{sec:lla}, this can be
rewritten as $D_0D_2$, a relation we have already discovered.
Thus six of the above degree two classes are independent.  Since the second
Betti number is six, these six classes generate the degree two part of 
$A^*(\mm_{0,2}(\PP^1,2))$, and so there are no other generators in degree two.

In degree three, any monomial with a factor of $D_1D_0$ can be expressed
in terms of other monomials.
Taking this into account together with the other
relations of lower degree, we have
six remaining degree three monomials in the divisor classes.  A generic 
relation among these has the form
\[aD_1^3+bD_1^2H_1+cD_1^2H_2+dD_1H_1H_2+eD_0^3+fD_0^2H_1=0\text{.}\]
Multiplying this expression by the four independent divisor classes and
integrating, we get a system of four linear equations in six variables.
The set of solutions is two-dimensional: any relation must have the form
\[aD_1^3+bD_1^2H_1+bD_1^2H_2+(-6a-4b)D_1H_1H_2+(32a-8b)D_0^3+(-96a-8b)D_0^2H_1
=0\text{.}\]
One can check that the relations
\begin{equation}
\label{eqn:c1}
(D_1+D_2)^3=2^3(D_1^3-3D_1^2H_1-3D_1^2H_2+6D_1H_1H_2+56D_0^3-72D_0^2H_1)
\end{equation}
and
\begin{equation}
\label{eqn:c2} 
D_1\psi_1\psi_2=\frac{1}{4}D_1^3-D_1^2H_1-D_1^2H_2+
\frac{5}{2}D_1H_1H_2+16D_0^3-16D_0^2H_1
\end{equation}
found above satisfy these conditions. Moreover, they
are clearly independent.
So they span the space of relations.
Thus four of these six degree three classes must be independent.  Since the
third Betti number is four, they generate the degree three part, and there
are no additional generators.

The degree four part is one-dimensional, so since $D_1^2H_1H_2$ is nonzero,
it generates the degree four part.

\subsection{A presentation for $A^*(\mm_{0,2}(\PP^1,2))$}

We have established the following result.

\begin{thm}\label{thm:prez}
With notation as established in Section \ref{sec:gen}, we have an isomorphism

\begin{equation}\label{geomprez}
A^*(\mm_{0,2}(\PP^1,2))\iso\frac{\QQ[D_0,D_1,D_2,H_1,H_2,\psi_1,\psi_2]}
{\left(\begin{array}{c}
H_1^2, H_2^2,D_0\psi_1,D_0\psi_2,D_2-\psi_1-\psi_2, \\
\psi_1-\frac{1}{4}D_1-\frac{1}{4}D_2-D_0+H_1, (D_1+D_2)^3, \\
\psi_2-\frac{1}{4}D_1-\frac{1}{4}D_2-D_0+H_2, D_1\psi_1\psi_2
\end{array}\right)}
\end{equation}
of graded rings.
\end{thm}

It is also possible, of course to give a presentation that avoids use of the
$\psi$-classes. See \cite{C4} for an example of such a presentation.
This example is more
efficient than (\ref{geomprez}) in the sense that it has the minimum possible number of 
generators (four) and relations (six). However, 
it is not very geometric; one would be hard-pressed to give a geometric
explanation for some of the relations. A goal of efficiency also leads
to some complicated relations, and this type of presentation will be 
difficult to generalize.  Including the $\psi$-classes as generators
leads to the geometric presentation (\ref{geomprez}), 
which is more beautiful and will
also be more useful.

\section{Computation of gravitational correlators using the presentation}
\label{sec:app}

Let $X$ be a smooth variety, $\b\in H_2(X)$, and $\g_1,\ldots,\g_n\in H^*(X)$.
The gravitational correlator 
$\bra\t_{d_1}\g_1,\ldots,\t_{d_n}\g_n\ket$ on $A^*(\mmm)$ is a rational number given by
\[\bra\t_{d_1}\g_1,\ldots,\t_{d_n}\g_n\ket_{g,\b}
=\int_{[\mmm]^\text{vir}}\prod_{i=1}^n
\left(\psi_i^{d_i}\ev_i^*(\g_i)\right)\text{.}\]
If $\b=0$, we must require either $g>0$ or $n\geq 3$. 
We usually suppress $\t_0$ from the notation, and we suppress the fundamental
class of $X$ or write it as $1$. We define any gravitational
correlator including an argument $\t_{-1}$ to be zero. Gravitational correlators
are important to physicists in the study of mirror symmetry.
Gromov-Witten invariants, for instance, are just
gravitational correlators with $d_i=0$ for all $i\in\und{n}$, so that there
are no $\psi$-classes in the corresponding integral.
The above information and much more can be found in \cite[Chapter 10]{CK}.

In this section, we will use the presentation given in Theorem \ref{thm:prez}
to compute all the genus zero, degree two, two-point gravitational correlators
of $\PP^1$. Algorithms for computing gravitational correlators have already
been constructed using indirect methods. We will show that our results agree
with the numbers computed by these existing methods. This provides a check
on the validity of the presentation.

Gravitational correlators are known to satisfy certain axioms, including the Degree,
Equivariance, Fundamental Class, Divisor, Splitting, and Dilaton Axioms. The 
algorithms mentioned above make use of these axioms in computing the correlators.
The Equivariance Axiom states that, if the $\g_i$ are homogeneous and $i\in\und{n-1}$, then
\begin{eqnarray*}
& & \bra\t_{d_1}\g_1,\ldots,\t_{d_{i+1}}\g_{i+1},
\t_{d_i}\g_{i},\ldots,\t_{d_n}\g_n\ket_{g,\b}\\
& = & (-1)^{\deg\g_i\cdot\deg\g_{i+1}}\bra\t_{d_1}\g_1,
\ldots,\t_{d_i}\g_{i},\t_{d_{i+1}}\g_{i+1},
\ldots,
\t_{d_n}\g_n\ket_{g,\b}\text{.}
\end{eqnarray*}
This has an obvious extension to any permutation of the entries, so that
``equivariance" refers to $S_n$-equivariance. Again, we consider only
cases where the cohomology lives
in even degrees, so that the gravitational correlators are in fact
{\em invariant} under permutation of the entries. We refer the reader to 
\cite[Chapter 10]{CK} for descriptions of the other axioms.

There are sixteen genus zero, degree two, two-point gravitational correlators
of $\PP^1$. We will compute them using the presentation given in Section
\ref{sec:prez}, invoking the Equivariance Axiom in
order to reduce the number of calculations to nine.
The computations were carried out using the algebraic geometry
software system Macaulay 2 (\cite{GS}) by entering the presentation for the
Chow ring and instructing the program to perform multiplications in this ring.
A sample of the code and its corresponding output can be found in \cite{C4}.

The top codimension ({\em i.e.} degree) 
piece of the Chow ring is generated by any non-trivial class
of that codimension.
Thus, once we know the degree of one such class, all integrals can be computed
in terms of it. The degrees of many such classes are given in Table
\ref{deg4ints}.
For our computations in Macaulay 2, it was convenient to use
the value
\[\int_{\mm}D_1^4=-20\text{.}\]
One may wish to avoid dependence on localization by using instead a
``geometrically obvious" degree. The value
\[\int_{\mm}D_2D_1H_1H_2=2\]
is appropriate for such a purpose. Indeed, it is rather clear that there
are two points in the moduli space that satisfy the conditions these
classes impose. They correspond to the following stable maps.

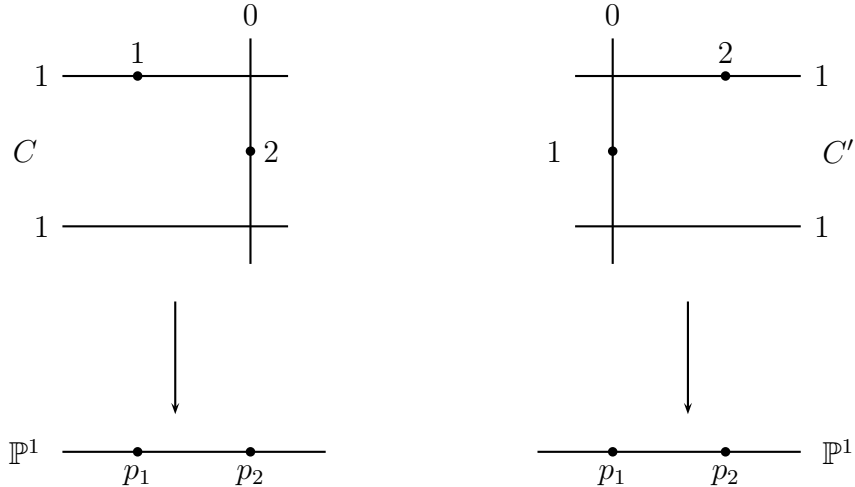
\begin{center}
\begin{pspicture}(0,-2.5)(4,4)
\rput(0,2){$C$}
\pnode(0.5,1){a}
\pnode(3.5,1){b}
\pnode(3,0.5){d}
\pnode(3,3.5){c}
\pnode(.5,3){f}
\pnode(3.5,3){e}
\dotnode(1.5,3){z}
\dotnode(3,2){y}
\ncline{a}{b}
\ncline{c}{d}
\ncline{e}{f}
\uput{5pt}[l](.5,1){1}
\uput{5pt}[l](.5,3){1}
\uput{5pt}[u](3,3.5){0}
\uput{5pt}[u](1.5,3){1}
\uput{5pt}[r](3,2){2}
\rput(0,-2){$\PP^1$}
\pnode(0.5,-2){g}
\pnode(4,-2){h}
\dotnode(1.5,-2){x}
\dotnode(3,-2){w}
\ncline{g}{h}
\uput{5pt}[d](3,-2){$p_2$}
\uput{5pt}[d](1.5,-2){$p_1$}
\pnode(2,0){i}
\pnode(2,-1.5){j}
\ncline{->}{i}{j}
\end{pspicture}
\hspace{1in}
%separate
\begin{pspicture}(0,-2.5)(4,4)
\rput(4,2){$C\pr$}
\pnode(0.5,1){a}
\pnode(3.5,1){b}
\pnode(1,0.5){d}
\pnode(1,3.5){c}
\pnode(.5,3){f}
\pnode(3.5,3){e}
\dotnode(2.5,3){z}
\dotnode(1,2){y}
\ncline{a}{b}
\ncline{c}{d}
\ncline{e}{f}
\uput{5pt}[r](3.5,1){1}
\uput{5pt}[r](3.5,3){1}
\uput{5pt}[u](1,3.5){0}
\uput{5pt}[u](2.5,3){2}
\uput{5pt}[l](0.5,2){1}
\rput(4,-2){$\PP^1$}
\pnode(0,-2){g}
\pnode(3.5,-2){h}
\dotnode(1,-2){x}
\dotnode(2.5,-2){w}
\ncline{g}{h}
\uput{5pt}[d](2.5,-2){$p_2$}
\uput{5pt}[d](1,-2){$p_1$}
\pnode(2,0){i}
\pnode(2,-1.5){j}
\ncline{->}{i}{j}
\end{pspicture}
\end{center}
Note also that these stable maps have no non-trivial automorphisms.

\renewcommand{\arraystretch}{1.5}
\begin{table}
\begin{center}
\begin{tabular}{|cccclcr|}
\hline
$\bra\t_4,1\ket_{0,2}$ & = & $\bra 1,\t_4\ket_{0,2}$ & = & $-20\cdot\frac{3}{80}$ 
& = & $-\frac{3}{4}$\\
$\bra\t_3H,1\ket_{0,2}$ & = & $\bra 1,\t_3H\ket_{0,2}$ & = & $-20\cdot-\frac{1}{80}$ 
& = & $\frac{1}{4}$\\
$\bra\t_3,H\ket_{0,2}$ & = & $\bra H,\t_3\ket_{0,2}$ & = & $-20\cdot\frac{1}{16}$ 
& = & $-\frac{5}{4}$ \\
$\bra\t_3,\t_1\ket_{0,2}$ & = & $\bra\t_1,\t_3\ket_{0,2}$ & = & $-20\cdot-\frac{3}{80}$ 
& = & $\frac{3}{4}$\\
$\bra\t_2H,\t_1\ket_{0,2}$ & = & $\bra\t_1,\t_2H\ket_{0,2}$ & = & $-20\cdot\frac{1}{80}$ 
& = & $-\frac{1}{4}$\\
$\bra\t_2H,H\ket_{0,2}$ & = & $\bra H,\t_2H\ket_{0,2}$ & = & $-20\cdot-\frac{1}{40}$ 
& = & $\frac{1}{2}$\\
& & $\bra\t_2,\t_2\ket_{0,2}$ & = & $-20\cdot-\frac{1}{16}$ 
& = & $\frac{5}{4}$\\
$\bra\t_2,\t_1H\ket_{0,2}$ & = & $\bra\t_1H,\t_2\ket_{0,2}$ & = & $-20\cdot\frac{3}{80}$ 
& = & $-\frac{3}{4}$\\ 
& & $\bra\t_1H,\t_1H\ket_{0,2}$ & = & $-20\cdot-\frac{1}{40}$ & = & $\frac{1}{2}$\\
\hline
\end{tabular}
\end{center}
\caption{Gravitational correlators via the presentation for 
$A^*(\mm_{0,2}(\PP^1,2))$
\label{gravcor}}
\end{table}
\renewcommand{\arraystretch}{1}

We show details for one example, the gravitational correlator
$\bra\t_2H,\t_1\ket_{0,2}$. By definition,
\[\bra\t_2H,\t_1\ket_{0,2}=\int_\mm \psi_1^2 H_1 \psi_2\text{.}\]
Macaulay 2 reduces the integrand to $\frac{1}{80}D^4$. Since $D^4$ has
degree $-20$, $\bra\t_2H_1,\t_1\ket=-\frac{1}{4}$. Similar computations result
in the values found in Table \ref{gravcor}. 

The values of these gravitational correlators can also be computed via procedures
described in \cite[Chapter 10]{CK} and results of Kontsevich and Manin in \cite{KM2}.
These calculations are carried out in \cite{C4}.
All of the values in Table \ref{gravcor} agree with those
computed by these previously established methods.

%\bibliographystyle{plain}
%\bibliography{global}

\end{document}